\documentclass[12pt,oneside,english]{amsart}
\usepackage[T1]{fontenc}
\usepackage[latin1]{inputenc}
\usepackage{geometry}
\usepackage{setspace}
\usepackage{amssymb}

\makeatletter
 \theoremstyle{plain}
\newtheorem{thm}{Theorem}[section]
  \theoremstyle{plain}
  \newtheorem{prop}[thm]{Proposition}
  \theoremstyle{plain}
  \newtheorem{lem}[thm]{Lemma}
  \theoremstyle{plain}
  \newtheorem{cor}[thm]{Corollary}

\numberwithin{equation}{section}

\usepackage{babel}
\makeatother
\begin{document}

\title{limit theorems for free multiplicative convolutions }

\author{hari bercovici and jiun-chau wang}

\begin{abstract}
We determine the distributional behavior for products of free random
variables in a general infinitesimal triangular array. The main theorems
in this paper extend a result for measures supported on the positive
half-line in \cite{BP in Mulimit}, and provide a new limit theorem
for measures on the unit circle with nonzero first moment. 
\end{abstract}
\maketitle

\section{Introduction}

Given two probability measures $\mu$, $\nu$ on $\mathbb{R}_{+}=(0,+\infty)$,
we will denote by $\mu\circledast\nu$ their classical multiplicative
convolution, and by $\mu\boxtimes\nu$ their free multiplicative convolution.
Thus, $\mu\circledast\nu$ is the distribution of $XY$, where $X$
and $Y$ are classically independent positive random variables with
distributions $\mu$ and $\nu$, respectively. Analogously, $\mu\boxtimes\nu$
is the distribution of $X^{1/2}YX^{1/2}$, where $X$ and $Y$ are
freely independent positive random variables with distributions $\mu$
and $\nu$. A triangular array $\{\nu_{nk}:\, n\geq1,\,1\leq k\leq k_{n}\}$
of probability measures on $\mathbb{R}_{+}$ is said to be \emph{infinitesimal}
if \[
\lim_{n\rightarrow\infty}\max_{1\leq k\leq k_{n}}\nu_{nk}(\{ t\in\mathbb{R}_{+}:\,\left|t-1\right|\geq\varepsilon\})=0,\]
for every $\varepsilon>0$. Given such an array, one is interested
in the asymptotic behavior of the measures\[
\mu_{n}=\nu_{n1}\circledast\nu_{n2}\circledast\cdots\circledast\nu_{nk_{n}}\]
and \[
\nu_{n}=\nu_{n1}\boxtimes\nu_{n2}\boxtimes\cdots\boxtimes\nu_{nk_{n}}.\]
The case of $\mu_{n}$ is completely understood, and is reduced to
the theory of addition of independent random variables by a logarithmic
change of variables. However, the free case $\nu_{n}$ does not simply
reduce to the additive theory by this change of variables.

The problem was first addressed in \cite{BP in Mulimit}, where a
triangular array such that $\nu_{n1}=\nu_{n2}=\cdots=\nu_{nk_{n}}$
for all $n$ was considered. In this case, necessary and sufficient
conditions were found for the weak convergence of the measures $\nu_{n}$.
In particular, it was shown that the sequence $\nu_{n}$ converges
weakly if $\mu_{n}$ converges, but not conversely.

In this paper we will find necessary and sufficient conditions for
the weak convergence of $\nu_{n}$ without any further assumptions
on the infinitesimal array. We also prove analogous results for infinitesimal
triangular arrays on the unit circle $\mathbb{T}=\{ z\in\mathbb{C}:\,\left|z\right|=1\}$.
In both cases the possible limit of $\nu_{n}$ is $\boxtimes$-infinitely
divisible as shown in \cite{CG1,Ser & Ber in free hincin}. 

The additive version of our results were studied earlier. Thus, consider
an array $\{\mu_{nk}:\, n\geq1,\,1\leq k\leq k_{n}\}$ of probability
measures on $\mathbb{R}$. Infinitesimality in this case means that\[
\lim_{n\rightarrow\infty}\max_{1\leq k\leq k_{n}}\mu_{nk}(\{ t\in\mathbb{R}:\,\left|t\right|\geq\varepsilon\})=0,\qquad\varepsilon>0.\]
Denote by \[
\lambda_{n}=\mu_{n1}*\mu_{n2}*\cdots*\mu_{nk_{n}}\]
the classical additive convolutions, and by \[
\rho_{n}=\mu_{n1}\boxplus\mu_{n2}\boxplus\cdots\boxplus\mu_{nk_{n}}\]
the free additive convolutions of these measures. When $\mu_{n1}=\mu_{n2}=\cdots=\mu_{nk_{n}}$
for all $n$, it was shown in \cite{BP in stable} that $\lambda_{n}$
converges weakly if and only if $\rho_{n}$ converges weakly. This
result was extended to arbitrary infinitesimal arrays by Chistyakov
and G\"{o}tze in \cite{CG2}. These authors made heavy use of analytic
subordination (first proved for the free additive convolution in \cite{V in free entropy}
generically and in \cite{Biane} for the general case; cf. also \cite{V in coalgebra,V in Markovian},
\cite{Ser & Ber in semigroup} and \cite{CG1} for different approaches).
Our methods do not require analytic subordination and are close to
the original approach in \cite{BP in Mulimit}.

The remainder of this paper is organized as follows. In Section 2,
we describe the analytical apparatus necessary for the calculation
of free multiplicative convolutions. We also describe the $\boxtimes$-infinitely
divisible measures on $\mathbb{R}_{+}$ and $\mathbb{T}$ and some
useful approximation results. In Sections 3 we give the convergence
criteria for arrays on $\mathbb{R}_{+}$, and in Section 4, we prove
the analogous result for $\mathbb{T}$.

\section{Preliminaries}

The analogue of Fourier transform for multiplicative free convolutions
was discovered by Voiculescu \cite{V in multiplication} (see also
\cite{BV in LevyHincin,BV in unbdd}). Denote by $\mathcal{M}_{+}$
the collection of Borel probability measures defined on $\mathbb{R}_{+}$,
and by $\mathcal{M}_{\mathbb{T}}^{\times}$ Borel probability measures
$\nu$ supported on the circle $\mathbb{T}$ with nonzero first moment,
i.e. $\int_{\mathbb{T}}t\, d\nu(t)\neq0$. 

Given $\nu\in\mathcal{M}_{+}$, one defines the analytic function
$\psi_{\nu}$ by\[
\psi_{\nu}(z)=\int_{0}^{\infty}\frac{tz}{1-tz}\, d\nu(t),\qquad z\in\mathbb{C}\setminus(0,+\infty).\]
The function $\psi_{\nu}$ is univalent in the left half-plane $i\mathbb{C}^{+}$,
and $\psi_{\nu}(i\mathbb{C}^{+})$ is a region contained in the circle
with diameter $(-1,0)$; moreover, $\psi_{\nu}(i\mathbb{C}^{+})\cap(-\infty,0)=(-1,0)$.
Setting $\Omega_{\nu}=\psi_{\nu}(i\mathbb{C}^{+})$, one defines the
\emph{$S$-transform} of the measure $\nu$ to be\[
S_{\nu}(z)=\frac{1+z}{z}\psi_{\nu}^{-1}(z),\qquad z\in\Omega_{\nu}.\]
The remarkable property of the $S$-transform is that for $\mu,\nu\in\mathcal{M}_{+}$,
one has\[
S_{\mu\boxtimes\nu}(z)=S_{\mu}(z)S_{\nu}(z),\]
for every $z$ in a neighborhood of $(-1,0)$.

For $\nu\in\mathcal{M}_{\mathbb{T}}^{\times}$, the function $\psi_{\nu}$
is defined by the formula given above (with the integral calculated
over $\mathbb{T}$), but its domain of definition is now the open
unit disk $\mathbb{D}=\{ z\in\mathbb{C}:\,\left|z\right|<1\}$. The
function $\psi_{\nu}$ has an inverse in a neighborhood of zero since
$\psi_{\nu}^{\prime}(0)=\int_{\mathbb{T}}t\, d\nu(t)\neq0$. The corresponding
$S$-transform is defined in a neighborhood of zero. It is sometimes
convenient to use a variation of the $S$-transform:\[
\Sigma_{\nu}(z)=S_{\nu}\left(\frac{z}{1-z}\right).\]
If $\nu\in\mathcal{M}_{\mathbb{T}}^{\times}$, the function $\Sigma_{\nu}$
is also defined in a neighborhood of zero, and \[
\Sigma_{\mu\boxtimes\nu}(z)=\Sigma_{\mu}(z)\Sigma_{\nu}(z),\qquad\mu,\nu\in\mathcal{M}_{\mathbb{T}}^{\times},\]
for all $z$ in a neighborhood of zero where all functions involved
are defined.

The weak convergence of probability measures can be translated in
terms of their $S$-transforms. 

\begin{thm}
\cite{BV in LevyHincin,BV in unbdd}
\begin{enumerate}
\item Given $\{\nu_{n}\}_{n=1}^{\infty}\subset\mathcal{M}_{+}$, the sequence
$\{\nu_{n}\}_{n=1}^{\infty}$ converges weakly to a measure $\nu\in\mathcal{M}_{+}$
if and only if there exist two positive numbers $0<b<a<1$ such that
the disk $D$ with the diameter $(-a,-b)$ is contained in $\Omega_{\nu_{n}}$
for all $n$, and the sequence $\{ S_{\nu_{n}}\}_{n=1}^{\infty}$
converges uniformly on $D$ to a function $S$. 
\item Given $\{\nu_{n}\}_{n=1}^{\infty}\subset\mathcal{M}_{\mathbb{T}}^{\times}$,
the sequence $\{\nu_{n}\}_{n=1}^{\infty}$ converges weakly to a measure
$\nu\in\mathcal{M}_{\mathbb{T}}^{\times}$ if and only if there exists
a neighborhood of zero $K\subset\mathbb{D}$ such that for all $\Sigma_{\nu_{n}}$
are defined in $K$, and the sequence $\{\Sigma_{\nu_{n}}\}_{n=1}^{\infty}$
converges uniformly on $K$ to a function $\Sigma$.
\end{enumerate}
Moreover, if \textup{(1)} is satisfied then $S=S_{\nu}$, and if \textup{(2)}
is satisfied then $\Sigma=\Sigma_{\nu}$. 

\end{thm}
An array $\{\nu_{nk}\}_{n,k}\subset\mathcal{M}_{\mathbb{T}}^{\times}$
is infinitesimal if\[
\lim_{n\rightarrow\infty}\max_{1\leq k\leq k_{n}}\nu_{nk}(\{ t\in\mathbb{T}:\,\left|\arg t\right|\geq\varepsilon\})=0,\]
for every $\varepsilon>0$; here the principal value of the argument
is used. The following proposition gives an approximation of the $S$-transform
(see Theorem 3.1 in \cite{BP in Mulimit} and Theorem 1.1 and Theorem
2.1 in \cite{Ser & Ber in free hincin}).

\begin{prop}
For $0<b<a<1$ and $\varepsilon\in(0,1)$, define $\overline{D}$
to be the closed disk with diameter $[-a,-b]$, and set $K_{\varepsilon}=\{ z\in\mathbb{C}:\,\left|z\right|\leq\varepsilon\}$. 
\begin{enumerate}
\item If an array $\{\nu_{nk}\}_{n,k}\subset\mathcal{M}_{+}$ is infinitesimal,
then the functions $S_{\nu_{nk}}$ are defined in $\overline{D}$
for sufficiently large $n$, and we have\[
S_{\nu_{nk}}(z)=1+\left[\int_{0}^{\infty}\frac{1-t}{1+z-tz}\, d\nu_{nk}(t)\right](1+u_{nk}(z)),\]
for all $z\in\overline{D}$, where $\lim_{n\rightarrow\infty}\max_{1\leq k\leq k_{n}}\left|u_{nk}(z)\right|=0$
uniformly on $\overline{D}$. 
\item If an array $\{\nu_{nk}\}_{n,k}\subset\mathcal{M}_{\mathbb{T}}^{\times}$
is infinitesimal, then $S_{\nu_{nk}}$ are defined in $K_{\varepsilon}$
when $n$ is large, and we have\[
S_{\nu_{nk}}(z)=1+\left[\int_{\mathbb{T}}\frac{1-t}{1+z-tz}\, d\nu_{nk}(t)\right](1+v_{nk}(z)),\]
for all $z\in K_{\varepsilon}$, where $\lim_{n\rightarrow\infty}\max_{1\leq k\leq k_{n}}\left|v_{nk}(z)\right|=0$
uniformly on $K_{\epsilon}$. 
\end{enumerate}
\end{prop}
A measure $\nu\in\mathcal{M}_{+}$ is said to be \emph{$\circledast$-infinitely
divisible} if, for each $n\in\mathbb{N}$, there exists a measure
$\nu_{n}\in\mathcal{M}_{+}$ such that\[
\nu=\underbrace{\nu_{n}\circledast\nu_{n}\circledast\cdots\circledast\nu_{n}}_{n\,\text{times}}.\]
The notion of \emph{$\boxtimes$-infinite divisibility} is defined
analogously. The study of $\circledast$-infinitely divisible measures
on $\mathbb{R}_{+}$ reduces (by a change of variable) to the study
of the usual $*$-infinitely divisible measures on $\mathbb{R}$.
The Fourier transform needs to be replaced by the \emph{Mellin-Fourier}
transform of a measure $\nu\in\mathcal{M}_{+}$ defined by \[
\Phi_{\nu}(s)=\int_{0}^{\infty}t^{is}\, d\nu(t),\qquad s\in\mathbb{R}.\]
The fundamental property of the Mellin-Fourier transform is that\[
\Phi_{\mu\circledast\nu}(s)=\Phi_{\mu}(s)\Phi_{\nu}(s).\]
A $\circledast$-infinitely divisible measure $\nu\in\mathcal{M}_{+}$
has the Mellin-Fourier transform\[
\Phi_{\nu}(s)=\exp\left[i\lambda s+\int_{0}^{\infty}\left(t^{-is}-1+\frac{is\log t}{\log^{2}t+1}\right)\frac{\log^{2}t+1}{\log^{2}t}\, d\rho(t)\right],\qquad s\in\mathbb{R},\]
where $\lambda\in\mathbb{R}$ and $\rho$ is a finite positive Borel
measure on $\mathbb{R}_{+}$. We use the notation $\nu_{\circledast}^{\lambda,\rho}$
to denote the $\circledast$-infinitely divisible measure determined
by $\lambda$ and $\rho$. For $\boxtimes$-infinite divisibility
we have the following formulas as in \cite{BV in LevyHincin,BV in unbdd}.
A measure $\nu\in\mathcal{M}_{+}$ is $\boxtimes$-infinitely divisible
if and only if there exist $\gamma\in\mathbb{R}$ and a finite positive
Borel measure $\sigma$ on the compact space $[0,+\infty]$ such that
$S_{\nu}(z)=\exp(v_{\gamma,\sigma}(z))$, where $v_{\gamma,\sigma}$
is given by \[
v_{\gamma,\sigma}\left(\frac{z}{1-z}\right)=\gamma-\sigma(\{+\infty\})z+\int_{[0,+\infty)}\frac{1+tz}{z-t}\, d\sigma(t),\qquad z\in\mathbb{C}\setminus[0,1].\]
A measure $\nu\in\mathcal{M}_{\mathbb{T}}^{\times}$ is $\boxtimes$-infinitely
divisible if and only if there exist $\gamma\in\mathbb{R}$ and a
finite positive Borel measure $\sigma$ on $\mathbb{T}$ such that
$\Sigma_{\nu}(z)=\exp(u_{\gamma,\sigma}(z))$, where $u_{\gamma,\sigma}$
is given by\[
u_{\gamma,\sigma}(z)=-i\gamma+\int_{\mathbb{T}}\frac{1+tz}{1-tz}\, d\sigma(t),\qquad z\in\mathbb{D}.\]
We denote by $\nu_{\boxtimes}^{\gamma,\sigma}$ the $\boxtimes$-infinitely
divisible measure determined by $\gamma$ and $\sigma$. There is
a unique $\boxtimes$-infinitely divisible measure $m$ on $\mathbb{T}$
such that its first moment is zero. This is the Haar, or normalized
arclength measure.

We conclude this section with a result which will be used repeatedly. 

\begin{lem}
Consider a sequence $\{ r_{n}\}_{n=1}^{\infty}\subset\mathbb{R}$
and two triangular arrays $\{ z_{nk}:\, n\geq1,\,1\leq k\leq k_{n}\}$,
$\{ w_{nk}:\, n\geq1,\,1\leq k\leq k_{n}\}$ of complex numbers. Assume
that
\begin{enumerate}
\item $\Im w_{nk}\geq0$, for $n\geq1$ and $1\leq k\leq k_{n}$. 
\item \[
z_{nk}=w_{nk}(1+\varepsilon_{nk}),\]
where\[
\varepsilon_{n}=\max_{1\leq k\leq k_{n}}\left|\varepsilon_{nk}\right|\]
converges to zero as $n\rightarrow\infty$.
\item There exists a positive constant $M$ such that for sufficiently large
$n$,\[
\left|\Re w_{nk}\right|\leq M\Im w_{nk}.\]

\end{enumerate}
Then the sequence $\{ r_{n}+\sum_{k=1}^{k_{n}}z_{nk}\}_{n=1}^{\infty}$
converges if and only if the sequence $\{ r_{n}+\sum_{k=1}^{k_{n}}w_{nk}\}_{n=1}^{\infty}$
converges. Moreover, the two sequences have the same limit.

\end{lem}
\begin{proof}
The assumptions on $\{ z_{nk}\}_{n,k}$ and $\{ w_{nk}\}_{n,k}$ imply
\begin{equation}
\left|\left(r_{n}+\sum_{k=1}^{k_{n}}z_{nk}\right)-\left(r_{n}+\sum_{k=1}^{k_{n}}w_{nk}\right)\right|\leq2(1+M)\varepsilon_{n}\left(\sum_{k=1}^{k_{n}}\Im w_{nk}\right),\label{eq:2.1}\end{equation}
and\begin{equation}
(1-\varepsilon_{n}-M\varepsilon_{n})\left(\sum_{k=1}^{k_{n}}\Im w_{nk}\right)\leq\left|\sum_{k=1}^{k_{n}}\Im z_{nk}\right|,\label{eq:2.2}\end{equation}
for sufficiently large $n$. If the sequence $\{ r_{n}+\sum_{k=1}^{k_{n}}z_{nk}\}_{n=1}^{\infty}$
converges to a complex number $z$, (2.2) implies that $\{\sum_{k=1}^{k_{n}}\Im w_{nk}\}_{n=1}^{\infty}$
is bounded, and then (2.1) shows that the sequence $\{ r_{n}+\sum_{k=1}^{k_{n}}w_{nk}\}_{n=1}^{\infty}$
also converges to $z$. Conversely, if $\{ r_{n}+\sum_{k=1}^{k_{n}}w_{nk}\}_{n=1}^{\infty}$
converges to $z$, then the sequence $\{\sum_{k=1}^{k_{n}}\Im w_{nk}\}_{n=1}^{\infty}$
is bounded and hence by (2.1) the sequence $\{ r_{n}+\sum_{k=1}^{k_{n}}w_{nk}\}_{n=1}^{\infty}$
converges to $z$ as well.
\end{proof}

\section{Free Multiplicative Convolution on $\mathbb{R}_{+}$}

Given an infinitesimal triangular array $\{\nu_{nk}:1\leq k\leq k_{n},\, n\in\mathbb{N}\}\subset\mathcal{M}_{+}$
and $\tau>0$, define positive numbers\[
b_{nk}=\exp\left(\int_{e^{-\tau}}^{e^{\tau}}\log t\, d\nu_{nk}(t)\right),\]
and measures $\nu_{nk}^{\circ}$ by\[
d\nu_{nk}^{\circ}(t)=d\nu_{nk}(b_{nk}t).\]
Obviously, $\max_{1\leq k\leq k_{n}}\left|b_{nk}-1\right|\rightarrow0$
as $n\rightarrow\infty$, and hence the array $\{\nu_{nk}^{\circ}\}_{n,k}$
is also infinitesimal. Define\[
g_{nk}(w)=\int_{0}^{\infty}\frac{t^{2}-1}{t^{2}+1}\, d\nu_{nk}^{\circ}\left(\frac{1}{t}\right)+\int_{0}^{\infty}\left[\frac{1+tw}{w-t}\right]\frac{(t-1)^{2}}{t^{2}+1}\, d\nu_{nk}^{\circ}\left(\frac{1}{t}\right),\]
for $w\in\mathbb{C}\setminus[0,+\infty)$. Note that $g_{nk}(\overline{w})=\overline{g_{nk}(w)}$
and $\Im g_{nk}(w)\leq0$ for all $w$ such that $\Im w>0$. 

\begin{lem}
For every compact set $K\subset\mathbb{C}^{+}\cap(i\mathbb{C}^{+})$
there exists a positive constant $M=M(\tau,K)$ such that for sufficiently
large $n$, we have\[
\left|\Re g_{nk}(w)\right|\leq M\left|\Im g_{nk}(w)\right|,\qquad w\in K,\,1\leq k\leq k_{n}.\]

\end{lem}
\begin{proof}
We assume for convenience that $\tau=1$. No generality is lost since
one can make a linear change of variable to modify the value of $\tau$.
By a change of variable we have\[
\int_{0}^{\infty}\frac{t^{2}-1}{t^{2}+1}\, d\nu_{nk}^{\circ}\left(\frac{1}{t}\right)=\int_{-\infty}^{\infty}\frac{1-e^{2x}}{1+e^{2x}}\, d\rho_{nk}(x+a_{nk})\]
and\[
\int_{0}^{\infty}\frac{(t-1)^{2}}{t^{2}+1}\, d\nu_{nk}^{\circ}\left(\frac{1}{t}\right)=\int_{-\infty}^{\infty}\frac{(1-e^{x})^{2}}{1+e^{2x}}\, d\rho_{nk}(x+a_{nk}),\]
where the probability measure $\rho_{nk}$ is defined as $d\rho_{nk}(x)=d\nu_{nk}(e^{x})$,
and $a_{nk}=\int_{\left|x\right|<1}x\, d\rho_{nk}(x)$. Note that
the family $\{\rho_{nk}\}_{n,k}$ is now an infinitesimal family of
probability measures on $\mathbb{R}$, and hence \[
\lim_{n\rightarrow\infty}\max_{1\leq k\leq k_{n}}\left|a_{nk}\right|=0.\]

We proceed by rewriting\begin{align*}
\int_{-\infty}^{\infty}\frac{1-e^{2x}}{1+e^{2x}}\, d\rho_{nk}(x+a_{nk})= & \int_{\left|x\right|<1}\left[\frac{1-e^{2(x-a_{nk})}}{1+e^{2(x-a_{nk})}}+(x-a_{nk})\right]\, d\rho_{nk}(x)\\
 & -\int_{\left|x\right|\geq1}a_{nk}\, d\rho_{nk}(x)+\int_{\left|x\right|\geq1}\frac{1-e^{2(x-a_{nk})}}{1+e^{2(x-a_{nk})}}\, d\rho_{nk}(x).\end{align*}
It is easy to see that\[
\left|1-e^{2(x-a_{nk})}+(x-a_{nk})+(x-a_{nk})e^{2(x-a_{nk})}\right|\leq60(x-a_{nk})^{2},\]
for $\left|x\right|<1$. Consequently, for all $n$ and $k$ we have\[
\left|\int_{\left|x\right|<1}\left[\frac{1-e^{2(x-a_{nk})}}{1+e^{2(x-a_{nk})}}+(x-a_{nk})\right]\, d\rho_{nk}(x)\right|\leq60\int_{-\infty}^{\infty}\frac{(1-e^{x})^{2}}{1+e^{2x}}\, d\rho_{nk}(x+a_{nk}).\]
 Since the family $\{\rho_{nk}\}_{n,k}$ is infinitesimal, there exists
$N\in\mathbb{N}$ such that\[
\left|a_{nk}\right|\leq\frac{1}{2},\]
 for all $n\geq N$, $1\leq k\leq k_{n}$. Note that \[
\frac{(1-e^{x})^{2}}{1+e^{2x}}\geq\frac{(1-\sqrt{e})^{2}}{1+e}\quad\text{ and}\quad5(1-e^{x})^{2}\geq\left|1-e^{2x}\right|,\]
for all $\left|x\right|\geq\frac{1}{2}$. We deduce that for $n\geq N$,
$1\leq k\leq k_{n}$ we have\begin{align*}
\left|\int_{\left|x\right|\geq1}a_{nk}\, d\rho_{nk}(x)\right| & \leq\int_{\left|x\right|\geq1}\, d\rho_{nk}(x)\\
 & \leq\frac{1+e}{(1-\sqrt{e})^{2}}\int_{-\infty}^{\infty}\frac{(1-e^{x})^{2}}{1+e^{2x}}\, d\rho_{nk}(x+a_{nk}),\end{align*}
and\begin{align*}
\left|\int_{\left|x\right|\geq1}\frac{1-e^{2(x-a_{nk})}}{1+e^{2(x-a_{nk})}}\, d\rho_{nk}(x)\right| & \leq5\int_{-\infty}^{\infty}\frac{(1-e^{x})^{2}}{1+e^{2x}}\, d\rho_{nk}(x+a_{nk}).\end{align*}
Therefore, for sufficiently large $n$,\[
\left|\int_{0}^{\infty}\frac{t^{2}-1}{t^{2}+1}\, d\nu_{nk}^{\circ}\left(\frac{1}{t}\right)\right|\leq74\int_{0}^{\infty}\frac{(t-1)^{2}}{t^{2}+1}\, d\nu_{nk}^{\circ}\left(\frac{1}{t}\right).\]
The compactness of the set $K$ implies the existence of positive
constants $M_{1}$ and $M_{2}$ such that\[
\left|\Re\left[\frac{1+tw}{w-t}\right]\right|\leq M_{1}\]
and \[
\left|\Im\left[\frac{1+tw}{w-t}\right]\right|=-\Im\left[\frac{1+tw}{w-t}\right]\geq M_{2},\]
for all $t\in(0,+\infty)$ and $w\in K$. Hence, we have for sufficiently
large $n$ and for $w\in K$, \begin{eqnarray*}
\left|\Re g_{nk}(w)\right| & = & \left|\int_{0}^{\infty}\frac{t^{2}-1}{t^{2}+1}\, d\nu_{nk}^{\circ}\left(\frac{1}{t}\right)+\int_{0}^{\infty}\Re\left[\frac{1+tw}{w-t}\right]\frac{(t-1)^{2}}{t^{2}+1}\, d\nu_{nk}^{\circ}\left(\frac{1}{t}\right)\right|\\
 & \leq & 74\int_{0}^{\infty}\frac{(t-1)^{2}}{t^{2}+1}\, d\nu_{nk}^{\circ}\left(\frac{1}{t}\right)+M_{1}\int_{0}^{\infty}\frac{(t-1)^{2}}{t^{2}+1}\, d\nu_{nk}^{\circ}\left(\frac{1}{t}\right)\\
 & \leq & \frac{(74+M_{1})}{M_{2}}\int_{0}^{\infty}M_{2}\frac{(t-1)^{2}}{t^{2}+1}\, d\nu_{nk}^{\circ}\left(\frac{1}{t}\right)\\
 & \leq & -\frac{(74+M_{1})}{M_{2}}\int_{0}^{\infty}\Im\left[\frac{1+tw}{w-t}\right]\frac{(t-1)^{2}}{t^{2}+1}\, d\nu_{nk}^{\circ}\left(\frac{1}{t}\right)\\
 & = & \frac{74+M_{1}}{M_{2}}\left|\Im g_{nk}(w)\right|.\end{eqnarray*}
The result follows with $M=(74+M_{1})/M_{2}$.
\end{proof}
Fix a closed disk $\overline{D}\subset i\mathbb{C}^{+}$ with diameter
$[-a,-b]$, where $0<b<a<1$. By Proposition 2.2, $S_{\nu_{nk}}$
is defined in $\overline{D}$ for large $n$. Setting $w=z/(1+z)$,
and using the identity\[
\frac{(w-1)(t-1)}{w-t}=\frac{t^{2}-1}{t^{2}+1}+\left[\frac{1+tw}{w-t}\right]\frac{(t-1)^{2}}{t^{2}+1},\]
we see that the function $S_{\nu_{nk}^{\circ}}$ admits the following
approximation: \[
S_{\nu_{nk}^{\circ}}\left(\frac{w}{1-w}\right)-1=g_{nk}(w)\left(1+u_{nk}\left(\frac{w}{1-w}\right)\right),\]
in another closed disk $\overline{D_{0}}=\{ z/(1+z):\, z\in\overline{D}\}$
with real center, and \[
\lim_{n\rightarrow\infty}\max_{1\leq k\leq k_{n}}\left|u_{nk}\left(\frac{w}{1-w}\right)\right|=0,\]
uniformly for all $w\in\overline{D_{0}}$. Note that \[
S_{\nu_{nk}^{\circ}}\left(\frac{w}{1-w}\right)=b_{nk}S_{\nu_{nk}}\left(\frac{w}{1-w}\right).\]
The infinitesimality of the array $\{\nu_{nk}\}_{n,k}$ also shows
that $S_{\nu_{nk}}(z)$ converges uniformly in $k$ and $z\in\overline{D}$
to $1$ as $n\rightarrow\infty$; indeed, $S_{\delta_{1}}\equiv1$.
Hence, for sufficiently large $n$, the principal branch of $\log S_{\nu_{nk}}(z)$
is defined in $\overline{D}$. Furthermore, since\[
\log w=w-1+o(\left|w-1\right|),\]
as $w\rightarrow1$, it is easy to see from Lemma 3.1 and Lemma 2.3
that we have the following result. Fix a real number $\gamma$, and
a finite positive Borel measure $\sigma$ on $[0,+\infty]$.

\begin{lem}
Let $\{\alpha_{n}\}_{n=1}^{\infty}$ be a sequence of positive real
numbers. Then the sequence of functions $\{-\log\alpha_{n}+\sum_{k=1}^{k_{n}}\log S_{\nu_{nk}}(z)\}_{n=1}^{\infty}$
converges to $v_{\gamma,\sigma}(z)$ uniformly for all $z\in\overline{D}$
as $n\rightarrow\infty$ if and only if \[
\lim_{n\rightarrow\infty}\left(-\log\alpha_{n}+\sum_{k=1}^{k_{n}}\left[g_{nk}(w)-\log b_{nk}\right]\right)=v_{\gamma,\sigma}\left(\frac{w}{1-w}\right)\]
uniformly for all $w\in\overline{D_{0}}$. 
\end{lem}
\begin{thm}
For an infinitesimal family $\{\nu_{nk}\}_{n,k}\subset\mathcal{M}_{+}$
and a sequence $\{\alpha_{n}\}_{n=1}^{\infty}\subset\mathbb{R}_{+}$
, the following two assertions are equivalent:
\begin{enumerate}
\item The sequence $\nu_{n1}\boxtimes\nu_{n2}\boxtimes\cdots\boxtimes\nu_{nk_{n}}\boxtimes\delta_{\alpha_{n}}$
converges weakly to $\nu_{\boxtimes}^{\gamma,\sigma}$. 
\item The sequence of measures\[
d\sigma_{n}(t)=\sum_{k=1}^{k_{n}}\frac{(t-1)^{2}}{t^{2}+1}\, d\nu_{nk}^{\circ}\left(\frac{1}{t}\right)\]
converges weakly in $[0,+\infty]$ to $\sigma$, and the sequence
\[
\gamma_{n}=-\log\alpha_{n}+\sum_{k=1}^{k_{n}}\left[\int_{0}^{\infty}\frac{t^{2}-1}{t^{2}+1}\, d\nu_{nk}^{\circ}\left(\frac{1}{t}\right)-\log b_{nk}\right]\]
converges to $\gamma$ as $n\rightarrow\infty$. 
\end{enumerate}
\end{thm}
\begin{proof}
Assume (1) holds. From Theorem 2.1, there exists a closed disk $\overline{D}$
with real center such that \[
\lim_{n\rightarrow\infty}\left(\frac{1}{\alpha_{n}}\prod_{k=1}^{k_{n}}S_{\nu_{nk}}(z)\right)=S_{\nu_{\boxtimes}^{\gamma,\sigma}}(z)=\exp(v_{\gamma,\sigma}(z))\]
uniformly on the disk $\overline{D}$. We may choose $\overline{D}$
small enough so that $\exp(v_{\gamma,\sigma}(z))$ is in $-i\mathbb{C}^{+}$
on $\overline{D}$. Applying the principal branch of the logarithm
function, we deduce that\[
\lim_{n\rightarrow\infty}\left(-\log\alpha_{n}+\sum_{k=1}^{k_{n}}\log S_{\nu_{nk}}(z)\right)=v_{\gamma,\sigma}(z),\]
uniformly on $\overline{D}$. Thus, Lemma 3.2 implies that\begin{equation}
\lim_{n\rightarrow\infty}\left(-\log\alpha_{n}+\sum_{k=1}^{k_{n}}\left[g_{nk}(w)-\log b_{nk}\right]\right)=v_{\gamma,\sigma}\left(\frac{w}{1-w}\right)\label{eq:3.1}\end{equation}
uniformly on $\overline{D_{0}}=\{ z/(1+z):\, z\in\overline{D}\}$.
Note that\begin{equation}
-\log\alpha_{n}+\sum_{k=1}^{k_{n}}\left[g_{nk}(w)-\log b_{nk}\right]=\gamma_{n}+\int_{0}^{\infty}\frac{1+tw}{w-t}\, d\sigma_{n}(t).\label{eq:3.2}\end{equation}
Considering the imaginary part of the equation (3.1), we have\begin{equation}
-\Im w\int_{[0,+\infty)}\frac{1+t^{2}}{\left|w-t\right|^{2}}\, d\sigma(t)=-\lim_{n\rightarrow\infty}\Im w\int_{(0,+\infty)}\frac{1+t^{2}}{\left|w-t\right|^{2}}\, d\sigma_{n}(t).\label{eq:3.3}\end{equation}
Note that the function $t\mapsto\frac{1+t^{2}}{\left|w-t\right|^{2}}$
is bounded away from zero and infinity for all $w\in\overline{D_{0}}$;
moreover, if $\Im w\ne0$ then (3.3) shows that\[
\sup_{n}\sigma_{n}((0,+\infty))<+\infty,\]
and hence the family $\{\sigma_{n}\}_{n=1}^{\infty}$ has a weak cluster
point $\sigma^{\prime}$ on the compact space $[0,+\infty]$. Then
(3.3) shows that $\sigma^{\prime}=\sigma$ , and consequently the
measures $\sigma_{n}$ converges weakly to $\sigma$ on $[0,+\infty]$
as $n\rightarrow\infty$. Then it is easy to see from (3.1) and (3.2)
that $\lim_{n\rightarrow\infty}\gamma_{n}=\gamma$.

Conversely, assume (2) holds. The infinitesimality of the array $\{\nu_{nk}\}_{n,k}$
implies that there exist $a^{\prime},b^{\prime}\in(0,1)$ with $b^{\prime}<a^{\prime}$
such that $S_{\nu_{nk}}$ are defined in $\overline{D^{\prime}}$,
the closed disk with the diameter $[-a^{\prime},-b^{\prime}]$, for
sufficiently large $n$. Let $\overline{D_{0}^{\prime}}=\{ z/(1+z):\, z\in\overline{D^{\prime}}\}$
and observe that there exists a positive constant $M=M(a^{\prime},b^{\prime})$
such that\[
\left|\frac{1+tw}{w-t}\right|\leq M,\qquad w\in\overline{D_{0}^{\prime}},\, t\in(0,+\infty).\]
 Thus, in view of (3.2), we deduce that (3.1) holds pointwise in $\overline{D_{0}^{\prime}}$.
Since $\Im g_{nk}(w)\leq0$ for $w\in D_{0}^{\prime}\cap i\mathbb{C}^{+}$,
the family $\{-\log\alpha_{n}+\sum_{k=1}^{k_{n}}\left[g_{nk}(w)-\log b_{nk}\right]\}_{n=1}^{\infty}$
is normal in $D_{0}^{\prime}\cap i\mathbb{C}^{+}$. Moreover, note
that $g_{nk}(\overline{w})=\overline{g_{nk}(w)}$ and\[
v_{\gamma,\sigma}\left(\frac{\overline{w}}{1-\overline{w}}\right)=\overline{v_{\gamma,\sigma}\left(\frac{w}{1-w}\right)},\]
for $w\in\overline{D_{0}^{\prime}}$. Therefore, as an application
of Montel's theorem, we conclude that (3.1) holds uniformly on compact
subsets of $\overline{D_{0}^{\prime}}$. From Lemma 3.2 we conclude
that there exists a smaller closed disk $\overline{D^{\prime\prime}}\subset\overline{D^{\prime}}$
with real center, such that\[
\lim_{n\rightarrow\infty}\left(-\log\alpha_{n}+\sum_{k=1}^{k_{n}}\log S_{\nu_{nk}}(z)\right)=v_{\gamma,\sigma}(z),\]
uniformly on $\overline{D^{\prime\prime}}$. Applying the exponential,
we obtain \[
\lim_{n\rightarrow\infty}\left(\frac{1}{\alpha_{n}}\prod_{k=1}^{k_{n}}S_{\nu_{nk}}(z)\right)=S_{\nu_{\boxtimes}^{\gamma,\sigma}}(z)=\exp(v_{\gamma,\sigma}(z))\]
uniformly on $\overline{D^{\prime\prime}}$. Therefore (1) follows
from Theorem 2.1.
\end{proof}
It has been pointed out in \cite{BP in Mulimit} that the weak convergence
criteria for products of free and independent random variables are
not equivalent. Nevertheless, the following correspondence is true.

\begin{cor}
Given an infinitesimal family $\{\nu_{nk}\}_{n,k}\subset\mathcal{M}_{+}$
and a sequence $\{\alpha_{n}\}_{n=1}^{\infty}\subset\mathbb{R}_{+}$,
the following two statements are equivalent:
\begin{enumerate}
\item The sequence $\nu_{n1}\boxtimes\nu_{n2}\boxtimes\cdots\boxtimes\nu_{nk_{n}}\boxtimes\delta_{\alpha_{n}}$
converges weakly to $\nu_{\boxtimes}^{\gamma,\sigma}$ and $\sigma(\{0\})=\sigma(\{\infty\})=0$;
\item The sequence $\nu_{n1}\circledast\nu_{n2}\circledast\cdots\circledast\nu_{nk_{n}}\circledast\delta_{\alpha_{n}}$
converges weakly to $\nu_{\circledast}^{\lambda,\rho}$. 
\end{enumerate}
If conditions \textup{(1)} and \textup{(2)} are satisfied then the
measure $\sigma$ and $\rho$ are related by\[
d\sigma(t)=\frac{\log^{2}t+1}{\log^{2}t}\frac{(t-1)^{2}}{t^{2}+1}\, d\rho(t),\]
and\[
\gamma-\lambda=\int_{0}^{\infty}\left(\frac{t^{2}-1}{t^{2}+1}+\frac{\log^{2}t}{\log^{2}t+1}\right)\frac{\log^{2}t+1}{\log^{2}t}\, d\rho(t).\]

\end{cor}
\begin{proof}
The proof is identical with that of Theorem 4.2 in \cite{BP in Mulimit}. 
\end{proof}

\section{Free Multiplicative Convolution on $\mathbb{T}$}

Fix an infinitesimal array $\{\nu_{nk}:\, n\geq1,\,1\leq k\leq k_{n}\}\subset\mathcal{M}_{\mathbb{T}}^{\times}$
and $\tau\in(0,\pi)$. Consider the centering constant\[
b_{nk}=\exp\left(\int_{\left|\arg t\right|<\tau}\log t\, d\nu_{nk}(t)\right),\]
and the centered measure $\nu_{nk}^{\circ}$ obtained as follows:\[
d\nu_{nk}^{\circ}(t)=d\nu_{nk}(b_{nk}t).\]
Here, as before, $\log t=i\arg t$ represents the principal branch
of $\log t$. We have $\max_{1\leq k\leq k_{n}}\left|\arg b_{nk}\right|\rightarrow0$
as $n\rightarrow\infty$, and hence the array $\{\nu_{nk}^{\circ}\}_{n,k}$
is also infinitesimal. Define \[
h_{nk}(z)=-i\int_{\mathbb{T}}\Im t\, d\nu_{nk}^{\circ}(t)+\int_{\mathbb{T}}\frac{1+tz}{1-tz}(1-\Re t)\, d\nu_{nk}^{\circ}(t),\qquad z\in\mathbb{D}.\]
Note that $\Re h_{nk}(z)\geq0$ for all $z\in\mathbb{D}$.

\begin{lem}
For every compact neighborhood of zero $K\subset\mathbb{D}$, there
exists a positive constant $M=M(\tau,K)$ such that for sufficiently
large $n$, we have \[
\left|\Im h_{nk}(z)\right|\leq M\Re h_{nk}(z),\qquad z\in K,\,1\leq k\leq k_{n}.\]

\end{lem}
\begin{proof}
We may again assume that $\tau=1$. Define probability measures $\rho_{nk}$
on $\mathbb{R}$ such that $\rho_{nk}(\sigma)=\nu_{nk}(e^{i\sigma})$
if $\sigma\subset[-\pi,\pi)$, and $\rho_{nk}(\sigma)=0$ when $\sigma\cap[-\pi,\pi)=\phi$.
Changing variables, we have\[
\int_{\mathbb{T}}\Im t\, d\nu_{nk}^{\circ}(t)=\int_{-\infty}^{\infty}\sin x\, d\rho_{nk}(x+a_{nk})\]
and\[
\int_{\mathbb{T}}(1-\Re t)\, d\nu_{nk}^{\circ}(t)=\int_{-\infty}^{\infty}(1-\cos x)\, d\rho_{nk}(x+a_{nk}),\]
where $a_{nk}=\int_{\left|x\right|<1}x\, d\rho_{nk}(x)=\int_{\left|\arg t\right|<1}\arg t\, d\nu_{nk}(t)$.
The infinitesimality of the family $\{\nu_{nk}\}_{n,k}$ implies that
for sufficiently large $n$, \[
\max_{1\leq k\leq k_{n}}\left|a_{nk}\right|\leq\frac{1}{10}.\]
By using the elementary inequalities\[
\left|\sin x-x\right|\leq2(1-\cos x),\qquad-2\leq x\leq2,\]
and \[
\frac{1}{10}+\left|\sin x\right|\leq10(1-\cos x),\,\text{where}\;\frac{9}{10}\leq\left|x\right|\leq\pi+\frac{9}{10},\]
we have, \begin{eqnarray*}
\left|\int_{\mathbb{T}}\Im t\, d\nu_{nk}^{\circ}(t)\right| & \leq & \left|\int_{\left|x\right|<1}\left[\sin(x-a_{nk})-(x-a_{nk})\right]\, d\rho_{nk}(x)\right|\\
 &  & +\left|\int_{[-\pi,-1]\cup[1,\pi)}a_{nk}\, d\rho_{nk}(x)\right|+\left|\int_{[-\pi,-1]\cup[1,\pi)}\sin(x-a_{nk})\, d\rho_{nk}(x)\right|\\
 & \leq & 12\int_{\mathbb{T}}(1-\Re t)\, d\nu_{nk}^{\circ}(t),\end{eqnarray*}
for large $n$ and $1\leq k\leq k_{n}$. Also, from the compactness
of the set $K$, there exist two positive constants $M_{1}$ and $M_{2}$
such that \[
\left|\Re\left[\frac{1+tz}{1-tz}\right]\right|=\Re\left[\frac{1+tz}{1-tz}\right]\geq M_{1}\]
and\[
\left|\Im\left[\frac{1+tz}{1-tz}\right]\right|\leq M_{2},\]
for all $t\in\mathbb{T}$ and $z\in K$. The result follows with $M=(12+M_{2})/M_{1}$. 
\end{proof}
Suppose $K\subset\mathbb{D}$ is a neighborhood of zero. The infinitesimality
of the array $\{\nu_{nk}\}_{n.k}$ implies that $S_{\nu_{nk}}(z)$
converges uniformly in $k$ and $z\in K$ to $1$ as $n\rightarrow\infty$,
and hence for sufficiently large $n$, $\Sigma_{\nu_{nk}}(z)$ and
the principal branch of $\log\Sigma_{\nu_{nk}}(z)$ are defined in
$K^{\prime}=\{ z/(1+z):\, z\in K\}$. 

Fix a real number $\gamma$, and a finite positive Borel measure $\sigma$
on $\mathbb{T}$. 

\begin{lem}
Let $\{\lambda_{n}\}_{n=1}^{\infty}\subset\mathbb{T}$. Then\[
\lim_{n\rightarrow\infty}\exp\left(-\log\lambda_{n}+\sum_{k=1}^{k_{n}}\log\Sigma_{\nu_{nk}}(z)\right)=\Sigma_{\nu_{\boxtimes}^{\gamma,\sigma}}(z)\]
uniformly on $K^{\prime}$ if, and only if \[
\lim_{n\rightarrow\infty}\exp\left(-\log\lambda_{n}+\sum_{k=1}^{k_{n}}\left[h_{nk}(z)-\log b_{nk}\right]\right)=\Sigma_{\nu_{\boxtimes}^{\gamma,\sigma}}(z)\]
uniformly on $K^{\prime}$. 
\end{lem}
\begin{proof}
From Proposition 2.2, we have the following approximation for the
function $S_{\nu_{nk}^{\circ}}$:\[
S_{\nu_{nk}^{\circ}}(z)=1+\left[\int_{\mathbb{T}}\frac{1-t}{1+z-tz}\, d\overline{\nu}_{nk}(t)\right](1+v_{nk}(z)),\qquad z\in K,\]
where\[
v_{n}(z)=\max_{1\leq k\leq k_{n}}\left|v_{nk}(z)\right|\]
satisfies $\lim_{n\rightarrow\infty}v_{n}(z)=0$ uniformly in $K$.
Introducing a change of variable $z\mapsto\frac{z}{1-z}$ and using
the identity\[
\frac{(1-t)(1-z)}{1-tz}=-i\Im t+\frac{1+tz}{1-tz}(1-\Re t),\]
we conclude that \[
b_{nk}\Sigma_{\nu_{nk}}(z)=\Sigma_{\nu_{nk}^{\circ}}(z)=1+h_{nk}(z)\left(1+v_{nk}\left(\frac{z}{1-z}\right)\right),\qquad z\in K^{\prime}.\]
Lemmas 4.1 and 2.3 imply that for any sequence of purely imaginary
numbers $\{ r_{n}\}_{n=1}^{\infty}$, the sequence $\{ r_{n}+\sum_{k=1}^{k_{n}}\left[\Sigma_{\nu_{nk}^{\circ}}(z)-1\right]\}_{n=1}^{\infty}$
converges uniformly on $K^{\prime}$ if and only if the sequence $\{ r_{n}+\sum_{k=1}^{k_{n}}h_{nk}(z)\}_{n=1}^{\infty}$
converges uniformly on $K^{\prime}$. Moreover, two sequences have
the same limit. Since $\log w/(w-1)\rightarrow1$ as $w\rightarrow1$,
we can replace $\Sigma_{\nu_{nk}^{\circ}}(z)-1$ by $\log\Sigma_{\nu_{nk}^{\circ}}(z)$.
The result follows by choosing $r_{n}=-\log\lambda_{n}-\sum_{k=1}^{k_{n}}\log b_{nk}$.
\end{proof}
\begin{thm}
For an infinitesimal array $\{\nu_{nk}\}_{n,k}\subset\mathcal{M}_{\mathbb{T}}^{\times}$
and a sequence $\{\lambda_{n}\}_{n=1}^{\infty}\subset\mathbb{T}$,
the following assertions are equivalent:
\begin{enumerate}
\item The sequence $\nu_{n1}\boxtimes\nu_{n2}\boxtimes\cdots\boxtimes\nu_{nk_{n}}\boxtimes\delta_{\lambda_{n}}$
converges weakly to $\nu_{\boxtimes}^{\gamma,\sigma}$.
\item The sequence of measures\[
d\sigma_{n}(t)=\sum_{k=1}^{k_{n}}(1-\Re t)\, d\nu_{nk}^{\circ}(t)\]
converges weakly on $\mathbb{T}$ to $\sigma$, and the limit\[
\lim_{n\rightarrow\infty}e^{i\gamma_{n}}=e^{i\gamma}\]
exists, where \[
\gamma_{n}=\arg\lambda_{n}+\sum_{k=1}^{k_{n}}\left[\int_{\mathbb{T}}\Im t\, d\nu_{nk}^{\circ}(t)+\arg b_{nk}\right].\]

\end{enumerate}
\end{thm}
\begin{proof}
Assume (1) holds. From Theorem 2.1, there exists $\varepsilon\in(0,1)$
such that all $\Sigma_{\nu_{nk}}$ are defined in $K_{\varepsilon}^{\prime}=\{ z/(1+z):\,\left|z\right|\leq\varepsilon\}$,
and we have\[
\lim_{n\rightarrow\infty}\left(\frac{1}{\lambda_{n}}\prod_{k=1}^{k_{n}}\Sigma_{\nu_{nk}}(z)\right)=\Sigma_{\nu_{\boxtimes}^{\gamma,\sigma}}(z)=e^{u_{\gamma,\sigma}(z)}\]
uniformly on $K_{\varepsilon}^{\prime}$. Hence, by Lemma 4.2 and
the definition of $u_{\gamma,\sigma}(z)$, we have\begin{gather}
\lim_{n\rightarrow\infty}\exp\left(-\log\lambda_{n}+\sum_{k=1}^{k_{n}}\left[h_{nk}(z)-\log b_{nk}\right]\right)=\exp\left(-i\gamma+\int_{\mathbb{T}}\left[\frac{1+tz}{1-tz}\right]\, d\sigma(t)\right)\label{eq:4.1}\end{gather}
 uniformly on $K_{\varepsilon}^{\prime}$. Taking the absolute value
on both sides of (4.1), we deduce that \begin{equation}
\lim_{n\rightarrow\infty}\exp\left(\Re\left[\sum_{k=1}^{k_{n}}h_{nk}(z)\right]\right)=\exp\left(\int_{\mathbb{T}}\Re\left[\frac{1+tz}{1-tz}\right]\, d\sigma(t)\right)\label{eq:4.2}\end{equation}
uniformly on $K_{\varepsilon}^{\prime}$. Note that\begin{equation}
-\log\lambda_{n}+\sum_{k=1}^{k_{n}}\left[h_{nk}(z)-\log b_{nk}\right]=-i\gamma_{n}+\int_{\mathbb{T}}\left[\frac{1+tz}{1-tz}\right]\, d\sigma_{n}(t).\label{eq:4.3}\end{equation}
Moreover, the real part of the function $-\log\lambda_{n}+\sum_{k=1}^{k_{n}}\left[h_{nk}(z)-\log b_{nk}\right]$
is the Poisson integral of the measure $d\sigma_{n}\left(\frac{1}{t}\right)$
and hence (4.2) uniquely determines the weak cluster point of $\{\sigma_{n}\}_{n=1}^{\infty}$
which is $\sigma$. We therefore conclude the weak convergence of
the sequence $\{\sigma_{n}\}_{n=1}^{\infty}$. Moreover, consider
$z=0$ in (4.1) and (4.2) to deduce that\[
\lim_{n\rightarrow\infty}\frac{e^{i\gamma}}{e^{i\gamma_{n}}}=1,\]
as desired. 

The converse implication is fairly easy now, since one can basically
reverse the steps to reach the statement (1) by using Lemma 4.2 and
the fact that $\{-\log\lambda_{n}+\sum_{k=1}^{k_{n}}\left[h_{nk}(z)-\log b_{nk}\right]\}_{n=1}^{\infty}$
is normal in $\mathbb{D}$. Therefore the details of the proof are
left to the reader.
\end{proof}
The previous result does not cover the possibility that the measures
$\nu_{n1}\boxtimes\nu_{n2}\boxtimes\cdots\boxtimes\nu_{nk_{n}}\boxtimes\delta_{\lambda_{n}}$
might converge to Haar measure $m$. We address now this special case.
Let us also note for further use the equality \[
\Sigma_{\nu}(0)=\frac{1}{\int_{\mathbb{T}}t\, d\nu(t)},\qquad\nu\in\mathcal{M}_{\mathbb{T}}^{\times}.\]

\begin{thm}
For an infinitesimal array $\{\nu_{nk}\}_{n,k}\subset\mathcal{M}_{\mathbb{T}}^{\times}$
and a sequence $\{\lambda_{n}\}_{n=1}^{\infty}\subset\mathbb{T}$,
the following assertions are equivalent:
\begin{enumerate}
\item The sequence $\nu_{n1}\boxtimes\nu_{n2}\boxtimes\cdots\boxtimes\nu_{nk_{n}}\boxtimes\delta_{\lambda_{n}}$
converges weakly to $m$.
\item \[
\lim_{n\rightarrow\infty}\sum_{k=1}^{k_{n}}\int_{\mathbb{T}}(1-\Re t)\, d\nu_{nk}^{\circ}(t)=+\infty.\]

\end{enumerate}
\end{thm}
\begin{proof}
Assume (2) holds. Define\[
\nu_{n}=\nu_{n1}\boxtimes\nu_{n2}\boxtimes\cdots\boxtimes\nu_{nk_{n}}\boxtimes\delta_{\lambda_{n}},\qquad n\in\mathbb{N}.\]
The compactness of $\mathbb{T}$ implies that $\{\nu_{n}\}_{n=1}^{\infty}$
is tight. Suppose $\nu$ is a weak cluster point of $\{\nu_{n}\}_{n=1}^{\infty}$.
From the free multiplicative analogue of Hin\v{c}in's theorem (Theorem
2.1 in \cite{Ser & Ber in free hincin}), the measure $\nu$ is $\boxtimes$-infinitely
divisible. By passing to a subsequence, we may assume that $\nu_{n}$
converges weakly to $\nu$ as $n\rightarrow\infty$. By (4.3), we
can reformulate the statement (2) as follows: \[
\lim_{n\rightarrow\infty}\Re\sum_{k=1}^{k_{n}}h_{nk}(0)=+\infty.\]
Then the inequality (2.2) implies that \[
\lim_{n\rightarrow\infty}\Re\sum_{k=1}^{k_{n}}\log\Sigma_{\nu_{nk}}(0)=+\infty,\]
and consequently we deduce that\begin{eqnarray*}
\left|\int_{\mathbb{T}}t\, d\nu(t)\right| & = & \lim_{n\rightarrow\infty}\left|\int_{\mathbb{T}}t\, d\nu_{n}(t)\right|\\
 & = & \lim_{n\rightarrow\infty}\frac{1}{\left|\prod_{k=1}^{k_{n}}\Sigma_{\nu_{nk}}(0)\right|}\\
 & = & \lim_{n\rightarrow\infty}\exp\left(-\Re\sum_{k=1}^{k_{n}}\log\Sigma_{\nu_{nk}}(0)\right)=0.\end{eqnarray*}
Therefore, the $\boxtimes$-infinitely divisible measure $\nu$ has
zero first moment, and hence we conclude that $\nu=m$. Moreover,
the full sequence $\nu_{n}$ converges weakly to $m$ since $\{\nu_{n}\}_{n=1}^{\infty}$
has a unique weak cluster point $m$.

Conversely, assume (1) holds but (2) fails to be true. By passing,
if necessary, to a subsequence, we may assume the sequence of measures
\[
d\sigma_{n}(t)=\sum_{k=1}^{k_{n}}(1-\Re t)\, d\nu_{nk}^{\circ}(t)\]
converges weakly to a finite positive Borel measure $\sigma$ on $\mathbb{T}$.
Since the sequence of functions $\{-\log\lambda_{n}+\sum_{k=1}^{k_{n}}\left[h_{nk}(z)-\log b_{nk}\right]\}_{n=1}^{\infty}$
is normal in $\mathbb{D}$, we may assume, by passing to a further
subsequence, that \[
\lim_{n\rightarrow\infty}\left(-\log\lambda_{n}+\sum_{k=1}^{k_{n}}\left[h_{nk}(z)-\log b_{nk}\right]\right)=f(z),\]
 uniformly on compact subsets of $\mathbb{D}$, where the function
$f$ is analytic in $\mathbb{D}$. Note that the function $f$ is
not identically infinity since \[
\Re f(z)=\lim_{n\rightarrow\infty}\int_{\mathbb{T}}\Re\left[\frac{1+tz}{1-tz}\right]\, d\sigma_{n}(t)=\int_{\mathbb{T}}\Re\left[\frac{1+tz}{1-tz}\right]\, d\sigma(t),\qquad z\in\mathbb{D}.\]
Then, as in the proof of Theorem 4.3, we conclude that there exists
$\gamma\in\mathbb{R}$ such that \[
\lim_{n\rightarrow\infty}e^{i\gamma_{n}}=e^{i\gamma},\]
where the number $\gamma_{n}$ is defined as in Theorem 4.3. An application
of Theorem 4.3 then shows that a subsequence of $\{\nu_{n}\}_{n=1}^{\infty}$
converges weakly to $\nu_{\boxtimes}^{\gamma,\sigma}$ which contradicts
(1). Therefore (2) must be true.
\end{proof}

\end{document}